\documentclass[12pt,oneside,english]{amsart}
\usepackage[latin1]{inputenc}
\usepackage{amssymb}

\makeatletter
\newtheorem{theorem}{Theorem}
\newtheorem{lemma}{Lemma}
\newtheorem{example}{Example}

\newtheorem{remark}{Remark}

\newtheorem{corollary}{Corollary}

\usepackage{babel}

\makeatother
\begin{document}
\baselineskip=17pt

\title[Binary additive problems: recursions ]{Binary additive problems: recursions for numbers of representations }

\author{Vladimir Shevelev}
\address{Department of Mathematics \\Ben-Gurion University of the
 Negev\\Beer-Sheva 84105, Israel. e-mail:shevelev@bgu.ac.il}

\subjclass{11P32}

\begin{abstract}
We prove some general recursions for the numbers of representations of positive integers as a sum
$x+y, x\in X, y\in Y,$ where $X,Y$ are increasing sequences. In particular, we obtain recursions for the number of  the Goldbach, Lemoine-Levy, Chen and other binary partitions.
\end{abstract}

\maketitle

\section{Introduction and main results}
Let $S=\{s_n\}_{n\geq1}, T=\{t_n\}_{n\geq1}$ be increasing sequences (finite or infinite) of positive odd numbers and $W=\{w_n\}_{n\geq1}$
 be increasing sequence of the common terms of sequences $S$ and $T.$ In this case we are writing $W=S\cap T.$ Denote $S(x), T(x),$  $W(x),$ etc., the corresponding counting functions of these sequences, e.g., $W(x)$ denotes $\#(n: w_n\leq x).$ For a positive even integer $x,$ denote $g(x)$ the number of solutions of the diophantine equation
\begin{equation}\label{1}
 x=s+t, \enskip s\in S, \;t\in T,
\end{equation}
when the order of the summands is not important; in particular, in the case of $s,t\in W$ we have only solution corresponding to $s,t.$
Everywhere below $s\in S,\enskip t\in T, \enskip w\in W ,$ etc.
In Section 2 we prove the following recursion.
\begin{theorem}\label{1}(Number of decompositions of even numbers into two odd numbers).\newline
For even $x\geq4,$
$$g(x)=\sum_{t_1\leq t\leq x/2}S(x-t)+\sum_{s_1\leq s\leq x/2}T(x-s)-\sum_{w_1\leq w\leq x/2}W(x-w)$$
\begin{equation}\label{2}
-S(x/2)T(x/2)+\begin{pmatrix} W(x/2)+1\\2\end{pmatrix}-g(x-2)-g(x-4)-...-g(2).
\end{equation}
\end{theorem}
\begin{corollary}\label{1}
If $S\subseteq T ,$ then, for even $x\geq4,$
$$g(x)=\sum_{t_1\leq t\leq x/2}S(x-t)+\sum_{s_1\leq s\leq x/2}(T(x-s)-S(x-s))$$\newpage
$$-S(x/2)T(x/2)+\begin{pmatrix} S(x/2)+1\\2\end{pmatrix}-g(x-2)-g(x-4)-...-g(2).$$
\end{corollary}
\begin{corollary}\label{2}
If $S=T,$ then, for even $x\geq4,$
$$g(x)=\sum_{s_1\leq s\leq x/2}S(x-s)
-\begin{pmatrix} S(x/2)\\2\end{pmatrix}-g(x-2)-g(x-4)-...-g(2).$$
\end{corollary}
Let, furthermore, $L=\{l_n\}_{n\geq1}, M=\{m_n\}_{n\geq1}$ be increasing sequences of positive even numbers and $Z=\{z_n\}_{n\geq1}=L\cap M$ and $L(x), M(x)$ and $Z(x)$ are their corresponding counting functions.For a positive even integer $x,$ denote $e(x)$ the number of solutions of the diophantine equation
$$ x=l+m, \enskip l\in L, m\in M,$$
when the order of the summands is not important.
\begin{theorem}\label{2}(Number of decompositions of even numbers into two even numbers).\newline
For even $x\geq2,$
$$e(x)=\sum_{m_1\leq m\leq x/2}L(x-m)+\sum_{l_1\leq l\leq x/2}M(x-l)-\sum_{z_1\leq z\leq x/2}Z(x-z)$$

$$-L(x/2)M(x/2)+\begin{pmatrix} Z(x/2)+1\\2\end{pmatrix}-e(x-2)-e(x-4)-...-e(0).$$

\end{theorem}
\begin{corollary}\label{3}
If $L\subseteq M ,$ then, for even $x\geq2,$
$$e(x)=\sum_{m_1\leq m\leq x/2}L(x-t)+\sum_{l_1\leq l\leq x/2}(M(x-s)-M(x-s))$$
$$-L(x/2)M(x/2)+\begin{pmatrix} L(x/2)+1\\2\end{pmatrix}-e(x-2)-e(x-4)-...-e(0).$$
\end{corollary}
\begin{corollary}\label{4}
If $L=M,$ then, for even $x\geq2,$
$$e(x)=\sum_{l_1\leq l\leq x/2}L(x-l)
-\begin{pmatrix} L(x/2)\\2\end{pmatrix}-e(x-2)-e(x-4)-...-e(0).$$
\end{corollary}
Finally, let  $U=\{u_n\}_{n\geq1} \enskip( V=\{v_n\}_{n\geq1})$  be increasing sequence of positive even (odd) numbers,
and $U(x) \enskip(V(x))$ be its counting function.
 \begin{theorem}\label{3}(Number of decompositions of odd numbers into one even and one odd numbers).\newline
 Let $h(x),$ for a positive odd integer $x,$  denote  the number of solutions of the diophantine equation \newpage
\begin{equation}\label{3}
 x=u+v, \enskip u\in U, v\in V.
\end{equation}
Then, for $x\geq3,$
$$h(x)=\sum_{v_1\leq v\leq (x+1)/2}U(x-v)+\sum_{u_1\leq u\leq (x+1)/2}V(x-u)$$
\begin{equation}\label{4}
-U((x+1)/2)V((x+1)/2)-h(x-2)-h(x-4)-...-h(1).
\end{equation}
\end{theorem}
In Section 2 we prove Theorems 1 - 3; in Section 3 we give some examples including recursions for numbers of the Goldbach, Lemoine-Levy, Chen and other binary partitions.
\section{Proof of theorems 1 - 3}
\bfseries Proof of Theorem 1.\enskip\mdseries a)We start with the base of induction of formula (2). Put in (2) $x=4.$ Let us prove that (2) is valid for every values of parameters $ s_1, t_1, w_1.$ For $x=4,$ denoting right hand side of (2)
via $r,$ we have
$$r=\sum_{t_1\leq t\leq 2}S(4-t)+\sum_{s_1\leq s\leq 2}T(4-s)-\sum_{w_1\leq w\leq 2}W(4-w)$$
\begin{equation}\label{5}
-S(2)T(2)+\begin{pmatrix} W(2)+1\\2\end{pmatrix}-g(2).
\end{equation}
We distinguish several cases:\newline
1) $s_1=1, t_1=1.$ Then also $w_1=1.$ Evidently, $g(2)=1.$ By (5), we have
$$ r=S(3)+T(3)-W(3)-1.$$
1a) $s_2=3, t_2=3,$ then $w_2=3.$  Therefore, $r=1.$ Evidently, we have only representation 4=3+1, and thus $r=g(4).$\newline
1b) $s_2=3, t_2\geq5,$ then $w_2\geq5.$ Again $r=1$ and $r=g(4).$\newline
1c) Symmetrically in case of $t_2=3, s_2\geq5$ we have $r=g(4)=1.$\newline
1d) $s_2\geq5, t_2\geq5,$ then $w_2\geq5$ and $r=0.$ Evidently, $g(4)=0$ and again $r=g(4).$\newline
2) $s_1=1, t_1=3,$ then $w_1\geq3$ and $g(2)=0.$ Now,according to (5), we have
$$r=T(3)=1=g(4).$$
3) Symmetrically in case of $t_1=1, s_1=3$ we have $r=S(3)=1=g(4).$\newline
4) Finally, if $s_1\geq3, t_1\geq3$ then $w\geq3$ and $g(2)=0$. Now , by (5), $r=0$ which corresponds to $g(4)=0.$ \newline
Thus the base of induction is valid.\newpage
b)Step of induction. Denote
$$ G(x)=\sum_{t_1\leq t\leq x/2}S(x-t)+\sum_{s_1\leq s\leq x/2}T(x-s)-\sum_{w_1\leq w\leq x/2}W(x-w)$$
\begin{equation}\label{6}
-S(x/2)T(x/2)+\begin{pmatrix} W(x/2)+1\\2\end{pmatrix}.
 \end{equation}
We should prove that
\begin{equation}\label{7}
G(x+2)-G(x)=g(x+2).
\end{equation}
We have
$$G(x+2)-G(x)=\sum_{t_1\leq t\leq x/2+1}(S(x+2-t)-S(x-t))$$
$$+\sum_{s_1\leq s\leq x/2+1}(T(x+2-s)-T(x-s))-\sum_{w_1\leq w\leq x/2+1}(W(x+2-w)-W(x-w))$$
$$+\begin{cases}
0,\;\;if\;\;x/2+1\not\in S \cup T,\\T(x/2-1),\;\;if\;\;x/2+1\in S\setminus W,\\S(x/2-1),\;\;if\;\;x/2+1\in T\setminus W,\\S(x/2-1)+T(x/2-1)-W(x/2-1),\;\;if\;\;x/2+1\in W\end{cases}$$
\begin{equation}\label{8}
-S(x/2+1)T(x/2+1)+S(x/2)T(x/2)+\begin{pmatrix} W(x/2+1)+1\\2\end{pmatrix}-\begin{pmatrix} W(x/2)+1\\2\end{pmatrix}.
\end{equation}
First of all, notice that if $x/2+1\in  S \cup T$ then $x/2$ is even and, consequently, $x/2\not\in  S \cup T.$
Therefore, we can rewrite (8) as follow:
$$G(x+2)-G(x)=\sum_{t_1\leq t\leq x/2+1}(S(x+2-t)-S(x-t))$$
$$+\sum_{s_1\leq s\leq x/2+1}(T(x+2-s)-T(x-s))-\sum_{w_1\leq w\leq x/2+1}(W(x+2-w)-W(x-w))$$
$$+\begin{cases}
0,\;\;if\;\;x/2+1\not\in S \cup T,\\T(x/2),\;\;if\;\;x/2+1\in S\setminus W,\\S(x/2),\;\;if\;\;x/2+1\in T\setminus W,\\S(x/2)+T(x/2)-W(x/2),\;\;if\;\;x/2+1\in W\end{cases}$$
\begin{equation}\label{9}
-S(x/2+1)T(x/2+1)+S(x/2)T(x/2)+\begin{pmatrix} W(x/2+1)+1\\2\end{pmatrix}-\begin{pmatrix} W(x/2)+1\\2\end{pmatrix}.
\end{equation}
Now we verify directly that
$$S(x/2+1)T(x/2+1)-S(x/2)T(x/2)=$$\newpage
\begin{equation}\label{10}
\begin{cases}
0,\;\;if\;\;x/2+1\not\in S \cup T,\\T(x/2),\;\;if\;\;x/2+1\in S\setminus W,\\S(x/2),\;\;if\;\;x/2+1\in T\setminus W,\\S(x/2)+T(x/2)+1,\;\;if\;\;x/2+1\in W\end{cases}
\end{equation}
and that
\begin{equation}\label{11}
\begin{pmatrix} W(x/2+1)+1\\2\end{pmatrix}-\begin{pmatrix} W(x/2)+1\\2\end{pmatrix}=W(x/2)+1.
\end{equation}
From (9)-(11) we find
$$G(x+2)-G(x)=\sum_{t_1\leq t\leq x/2+1}(S(x+2-t)-S(x-t))$$
$$+\sum_{s_1\leq s\leq x/2+1}(T(x+2-s)-T(x-s))-\sum_{w_1\leq w\leq x/2+1}(W(x+2-w)-W(x-w))=$$
\begin{equation}\label{12}
\sum_{t_1\leq t\leq x/2+1: \enskip x+2-t\in S}1+\sum_{s_1\leq s\leq x/2+1: \enskip x+2-s\in T}1-\sum_{w_1\leq w\leq x/2+1: \enskip x+2-w\in W}1.
\end{equation}
Notice that, if $s\leq x/2+1,$ then $t=x+2-s\geq x/2+1,$ if $t\leq x/2+1,$ then $s=x+2-t\geq x/2+1,$ and if
$s= x/2+1,$ then $t=x+2-s= x/2+1\in W, \enskip t= x/2+1,$ then $s=x+2-t= x/2+1\in W.$ Therefore, sum (12), indeed, is
$g(x+2). \blacksquare$ \newline
\bfseries Proof of Theorem 2\enskip\mdseries is quite analogous to proof of Theorem 1.\newline
\bfseries Proof of Theorem 3.\enskip\mdseries a)Put in (4) $x=3.$ Let us prove that (4) is valid for every values of parameters $ u_1, v_1.$ For $x=3,$ denoting right hand side of (4)
via $r,$ we have
\begin{equation}\label{13}
r=\sum_{v_1\leq v\leq 2}U(3-v)+\sum_{u_1\leq u\leq 2}V(3-u)-U(2)V(2)-h(1).
\end{equation}
We distinguish several cases:\newline
1)$u_1=0, v_1=1.$ Then $h(1)=1, V(2)=1$ and, by (13), we have
\begin{equation}\label{14}
r=V(3)-1+\begin{cases}1,\;\;if\;\;u_2=2,\\0\;\;if\;\;u_2\geq4.\end{cases}
\end{equation}
1a)$u_2=2, v_2=3,$ then $V(3)=2$ and by (14) we have $r=2.$ In this case we have two representations of 3: 3=0+3
and 3=2+1. Thus $r=h(3).$\newline
1b)$u_2=2, v_2\geq5,$ then $V(3)=1$ and, by (14),  $r=1.$ Here there is only representation: 3=2+1 and again $r=h(3).$\newline
1c)$u_2\geq4, v_2=3,$ then $V(3)=2$ and, by (14),  $r=1.$ Here there is only representation: 3=0+3 and again $r=h(3).$\newpage
1d)$u_2\geq4, v_2\geq5,$ then $V(3)=1$ and ,by (14),$r=0.$ Evidently, in this case $h(3)=0$ as well.\newline
2)$u_1=2, v_1=1.$ Then $h(1)=0, V(2)=U(2)=1$ and, by (13), we have $r=V(1)=1.$ Evidently, in this case $h(3)=1$ as well.\newline
3)$u_1=0, v_1=3.$ Then $h(1)=0, V(2)=U(2)=1$ and, by (13), we have $r=V(1)=1.$ Evidently, in this case $h(3)=1$ as well.\newline
4)$u_1\geq4, v_1\geq1$ or $u_1\geq0, v_1\geq5,$ then evidently $r=h(3)=0.$
b)Step of induction. Denote
$$ H(x)=\sum_{v_1\leq v\leq (x+1)/2}U(x-v)+\sum_{u_1\leq u\leq (x+1)/2}V(x-u)$$
\begin{equation}\label{15}
-U((x+1)/2)V((x+1)/2).
 \end{equation}
We should prove that
\begin{equation}\label{16}
H(x+2)-H(x)=h(x+2).
 \end{equation}
We have
$$H(x+2)-H(x)=\sum_{v_1\leq v\leq (x+1)/2}(U(x+2-v)-U(x-v))$$
$$+\sum_{u_1\leq u\leq (x+1)/2}(V(x+2-u)-V(x-u))$$
$$+\begin{cases}
0,\;\;if\;\;(x+3)/2\not\in U \cup V,\\U((x+1)/2),\;\;if\;\;(x+3)/2\in V,\\V(x+1)/2),\;\;if\;\;(x+3)/2\in U\end{cases}$$
\begin{equation}\label{17}
-U((x+3)/2)V((x+3)/2)+U((x+1)/2)V((x+1)/2).
\end{equation}
Now, taking into account that $U\cap V={\o},$ we verify directly that
$$U((x+3)/2)V((x+3)/2)-U((x+1)/2)V((x+1)/2)=$$
$$+\begin{cases}
0,\;\;if\;\;(x+3)/2\not\in U \cup V,\\U((x+1)/2),\;\;if\;\;(x+3)/2\in V,\\V(x+1)/2),\;\;if\;\;(x+3)/2\in U\end{cases}$$
and from (17) we find
$$H(x+2)-H(x)=\sum_{v_1\leq v\leq (x+1)/2}(U(x+2-v)-U(x-v))$$
$$+\sum_{u_1\leq u\leq (x+1)/2}(V(x+2-u)-V(x-u))=$$
\begin{equation}\label{18}
\sum_{v_1\leq v\leq (x+1)/2: \enskip x+2-v\in U}1+\sum_{u_1\leq u\leq (x+1)/2: \enskip x+2-u\in V}1.
\end{equation}
\newpage
Notice that, if $v\leq (x+1)/2,$ then $u=x+2-v> (x+1)/2;$  if $u\leq (x+1)/2,$ then $v=x+2-u> (x+1)/2.$
Therefore, sum (18) is, indeed, $h(x+2).\blacksquare$
\section{Some examples}
\begin{example}\label{1}
Number of decompositions of \enskip $2n \enskip(n\geq1)$ into unordered sums of two odd primes (Goldbach partitions)
\end{example}
The first terms of this sequence are (see A002375 in [1]):
\begin{equation}\label{19}
0, 0, 1, 1, 2, 1, 2, 2, 2, 2, 3, 3, 3, 2, 3, 2, 4, 4, 2, 3, 4, 3, 4, 5, 4, 3, 5, 3, 4, 6, ...
\end{equation}
Let $\pi (x)$ be the counting prime function, $ \pi_1 (x)$ be the counting function of odd primes. Note that,
Hardy and Wright give a simple explicit formula for $\pi (n)$ for $n\geq4$ (see [2]):
$$ \pi (n)=-1+\sum_{j=3}^{n}((j-2)!-j\lfloor (j-2)!/j\rfloor).$$
Denote $P=\{p\}$ the set of all primes, $ P_1$ the set of odd primes. Putting in Corollary 2 $x=2n$ and $S=P_1,$
we have a recursion
$$g(2n)=\sum_{3\leq p\leq n}\pi_1(2n-p)
-\begin{pmatrix} \pi_1(n)\\2\end{pmatrix}$$
\begin{equation}\label{20}
-g(2n-2)-g(2n-4)-...-g(2),\enskip n\geq2, \enskip g(2)=0.
\end{equation}
\begin{example}\label{2}
Number of decompositions of \enskip $2n \enskip(n\geq1)$ into sums of a prime and a prime or semiprime (Chen partitions)
\end{example}
Let $P_2=\{p_1p_2\}$ be the set of semiprimes, i.e. products of two (possibly equal) primes. Denote $P_{2,1}=\{p_{2,1}\}$ the subset of odd semiprimes. If $\pi_2(x)$ is the counting function of $P_2,$ then, by Noel-Panos-Wilson formula (see [2]), we have
\begin{equation}\label{21}
\pi_2(x)=\sum_{p\leq \sqrt{x},\enskip p\in P}(\pi(x/p)-\pi(p)+1).
\end{equation}
Denoting $\pi_{2,1}(x)$ the counting function of set $P_{2,1},$ by (21) we have
\begin{equation}\label{22}
\pi_{2,1}(x)=\sum_{3\leq p\leq \sqrt{x},\enskip p\in P}(\pi(x/p)-\pi(p)+1).
\end{equation}
Putting in Corollary 1 $S=P_1, T=P_1\cup P_{2,1}, x=2n,$ we obtain a recursion for the numbers of Chen "odd-odd" partitions:\newpage
$$g_1(2n)=\sum_{3\leq t\leq n,\;t\in P_1\cup P_{2,1}}\pi_1(2n-t)+
\sum_{3\leq p\leq n,\; p\in P}\pi_{2,1}(2n-p)
$$

$$-\pi_1(n)(\pi_1(n)/2-1/2+\pi_{2,1}(n))$$
\begin{equation}\label{23}
-g_1(2n-2)-...-g_1(2), \enskip n\geq2, g_1(2)=0.
\end{equation}
The first terms of this sequence are:
\begin{equation}\label{24}
0, 0, 1, 1, 2, 2, 3, 3, 3, 4, 5, 4, 6, 6, 4, 6, 6, 6, 8, 7, 7, ...
\end{equation}
If to take into account that the number of Chen "even-even" partitions is evidently
$$g_2(2n)=\begin{cases}
1,\;\;if\;\;n-1 \enskip is \enskip prime \enskip or \enskip 1,\\0,\;\;
otherwise\end{cases}$$
which could be written by the recursion: $g_2(2)=0,$ for $n\geq2,$
$$g_2(2n)=\pi(n-1)+1-g_2(2n-2)-...-g_2(2), $$
then the first terms of all Chen partitions $g(2n)=g_1(2n)+g_2(2n)$ are the following
 (see A155216 in [1], where $a(n):=g(2n)$):
$$0, 1, 2, 2, 2, 3, 3, 4, 3, 4, 5, 5, 6, 7, 4, 6, 6, 7, 8, 8, 7,...$$
with the recursion: $g(2)=0,$ for $n\geq2,$
$$g(2n)=\sum_{3\leq t\leq n,\;t\in P_1\cup P_{2,1}}\pi_1(2n-t)+
\sum_{3\leq p\leq n,\; p\in P}\pi_{2,1}(2n-p)-$$

 $$\pi_1(n)(\pi_1(n)/2-1/2+
 \pi_{2,1}(n))+$$
\begin{equation}\label{25}
  \pi(n-1)+1-g(2n-2)-...-g(2).
  \end{equation}
  \begin{remark}\label{r1} One can prove also another interesting recursion for the number of
   Chen's partitions:
$$g(2)=0; \enskip for \; n\geq2,\; g(2n) = \sum_{3\leq p\leq n,\; p\in P}\pi_1(2n - p) + 
\sum_{9\leq t\leq2n,\; t\in P_{2,1}}\pi_1(2n - t) +$$
 $$ \pi_1(n) - \binom {\pi_1(n)}{2} + \delta(n) - g(2n-2) - ... - g(2),$$
  where  $$\delta(n) =\begin{cases}
1,\;\;if\;\;n\in P, \\2,\;\;
otherwise.\end{cases}$$ 
It could be proven by induction in the same way as the proof of Theorem 1.
   \end{remark}
   \newpage
\begin{example}\label{3}
Number of decompositions of \enskip $2n+1 \enskip(n\geq0)$ into sums of a prime and a doubled prime (Lemoine-Levy partitions)
\end{example}
Put in Theorem 3 $U=2P, V=P_1, x=2n+1.$ Then $U(x)=\pi(x/2), V(x)=\pi_{1}(x)$ and we
 find a recursion for required numbers:
$$h(2n+1)=\sum_{3\leq p\leq n+1,\enskip p \enskip prime}\pi((2n+1-p)/2)+\sum_{2\leq q\leq
(n+1)/2, \enskip q\enskip prime}\pi_1(2n+1-2q)-$$
\begin{equation}\label{26}
\pi((n+1)/2)\pi_{1}(n+1)-h(2n-3)-...-h(1),\enskip n\geq2, h(1)=0.
\end{equation}
The first terms of this sequence are the following (see A046927 in [1], where
 $a(n):=h(2n+1)$)
\begin{equation}\label{27}
 0, 0, 0, 1, 2, 2, 2, 2, 4, 2, 3, 3, 3, 4, 4, 2, 5, 3, 4, 4, 5, 4, 6, 4, 4, 7,...
\end{equation}

\begin{example}\label{4}
Number of decompositions of \enskip $4n+1 \enskip(n\geq0)$ into unordered sums of two squares of nonnegative integers and number of decompositions of \enskip $n \enskip(n\geq0)$ into unordered sums of two triangular numbers
\end{example}
Let $Q_1=\{q_1\}(Q_2=\{q_2\})$ be the set of odd (even) squares of integers. Then its counting function is $\lfloor(\sqrt{x}+1)/2\rfloor\enskip (\lfloor(\sqrt{x}+2)/2\rfloor).$ Putting in Theorem 3 $U=Q_2, V=Q_1, x=4n+1$ and taking into account that $h(4i-1)=0, \enskip i\geq1,$ we obtain the following recursion:
$$h(4n+1)=\sum_{0\leq q_2\leq 2n}\lfloor(1+\sqrt{4n-q_2+1})/2\rfloor+\sum_{1\leq q_1\leq 2n+1}\lfloor(2+\sqrt{4n-q_1+1})/2\rfloor$$
$$-\lfloor(1+\sqrt{2n+1})/2\rfloor\lfloor(2+\sqrt{2n+1})/2\rfloor-h(4n-3)-h(4n-7)-...-h(1),$$
\begin{equation}\label{28}
n\geq1, h(1)=1.
\end{equation}
With help of the following statement we shall give a simpler recursion.
\begin{lemma}\label{1}
For $n\geq0,$ the number $4n+1$ is a sum of two squares if and only if $n$ is a sum of two triangular numbers.
Moreover, the number of unordered representations of $4n+1$ as a sum of two squares of nonnegative integers
equals to the number of unordered representations of $n$ as a sum of two triangular numbers.
\end{lemma}
\bfseries Proof.\enskip\mdseries For integer $x\geq0,$ denote $T_x=x(x+1)/2.$ If $n=T_x+T_y,$ then one can verify
directly that
$$4n+1=(x+y+1)^2+(x-y)^2.$$
Conversely, if $4n+1=x^2+y^2,$ then one can verify that
$$n=T_{(x+y-1)/2}+T_{(|x-y|-1)/2}.\blacksquare$$
\newpage
\begin{remark}\label{r2} Notice that, the identity
$$4T_x+4T_y+1=(x+y+1)^2+(x-y)^2$$
is a generalization of well known identity $T_x+T_{x-1}=x^2.$
\end{remark}
Let now $L=\{l\}_{l\geq0},$ where $l=j(j+1),\enskip j\geq0.$ Note that, $L(x)=\lfloor (1+\sqrt{4x+1})/2\rfloor.$ Consider representations of the form $2n=l_1+l_2, \enskip l_1,l_2\in L,\enskip n\geq0$ or, equivalently, $n=T_x+T_y.$ The number $t(n)=e(2n)$ of such representations we find by Corollary 4 for $x=2n.$
After small transformation we obtain the following recursion:
$$t(n)=\sum_{k=1}^{\lfloor(1+\sqrt{4n+1})/2\rfloor}\lfloor(1+\sqrt{1+4(2n-k^2+k)})/2\rfloor-\begin{pmatrix} \lfloor(1+\sqrt{4n+1})/2\rfloor\\2\end{pmatrix}$$
\begin{equation}\label{29}
-t(n-1)-t(n-2)-...-t(0),\enskip n\geq1, \enskip t(0)=1.
\end{equation}
By Lemma1, $$h(4n+1)=t(n).$$
The first terms of this sequence are (see A052343 in [1]):
\begin{equation}\label{30}
 1, 1, 1, 1, 1, 0, 2, 1, 0, 1, 1, 1, 1, 1, 0, 1, 2, 0, 1, 0, 1, 2, 1, 0, 1, 1, ...
\end{equation}

\end{document}